\documentstyle[11pt,amssymb]{amsart}

\swapnumbers

\def \a{\alpha}
\def \be{\beta}
\def \ga{\gamma}

\def \de{\delta}
\def \e{\varepsilon}
\def \la{\lambda}

\def \vr{\varphi}

\def \re{{\Bbb R}}
\def \na{{\Bbb N}}
\def \Z {{\Bbb Z}}

\def \lv{\left\vert}
\def \rv{\right\vert}
\def \lV{\left\Vert}
\def \rV{\right\Vert}
\def \ov{\overline}

\def \SS{{\Sigma^+}}
\def \s{\sigma}
\def \supp{{\text{\rm supp}}}
\def \fS{{\mathfrak S}}
\def \Fix{{\text{\rm Fix}}}
\def \0{{\text{\rm\o}}}
\def \b{\big}
\def \hA{{\widehat{A}}}
\def \cA{{\mathcal{A}}}
\def \cD{{\mathcal{D}}}
\def \E{{\text{\large \rm e}}}
\def \cF{{\mathcal F}}
\def \cG{{\mathcal{G}}}
\def \cH{{\mathcal{H}}}
\def \cL{{\mathcal L}}
\def \cM{{\mathcal{M}}}

\def \cO{{\mathcal{O}}}
\def \cU{{\mathcal{U}}}
\def \cV{{\mathcal{V}}}

\def \cY{{\mathcal Y}}
\def \Can{{C^\a(\Sigma^+,\re)}}
\def \Ca{{C^\a_0(\Sigma^+,\re)}}
\def \Cb{{C^\be_0(\Sigma^+,\re)}}
\def \Cga{{C^\ga_0(\Sigma^+,\re)}}
\def \Cap{{C^{\a+}_0(\Sigma^+,\re)}}

\def \tnu{{\widetilde{\nu}}}
\def \hmu{{\widehat{\mu}}}
\def \ta{{\theta_A}}
\def \td{{\widetilde{d}}}
\def \tf{{\theta_f}}
\def \tg{{\theta_g}}
\def \Hl{{{\mathcal H}_\lambda}}
\def \Hlo{{{\mathcal H}_\la^0}}
\def \hm{{\widehat{\mu}}}
\def \diam{{\text{\rm diam}}}
\def \zum{\textstyle\sum\limits}
\def \Hold{{\text{\rm Hold}}}

\newcommand{\abs}[1]{\left\vert#1\right\vert}
\newcommand{\norm}[1]{\left\Vert#1\right\Vert}

\theoremstyle{plain}
\newtheorem{Thm}{\bf Theorem}[section]

\newtheorem{Lemma}[Thm]{\bf Lemma}

\newtheorem{Corollary}[Thm]{\bf Corollary}
\newtheorem{Theorem}[Thm]{\bf Theorem}
\newtheorem{Proposition}[Thm]{\bf Proposition}

\theoremstyle{definition}
\newtheorem{Remark}[Thm]{\bf Remark}

\title{Maximizing measures for expanding transformations.}

\date{February, 1998.}

\author[G. Contreras]{Gonzalo Contreras}
\address{CIMAT \\
         P.O.Box 402\\
         36.000 Guanajuato, Gto.\\
         M\'exico.}
\email{gonzalo@@fractal.cimat.mx}         
\thanks{G. Contreras and A. Lopes were partially supported by
        CNPq-Brazil.} 

\author[A. Lopes]{Artur Lopes}
\address{Instituto de Matem\'atica \\
         UFRGS \\
         Porto Alegre. Brasil.}
\email{alopes@@mat.ufrgs.br}         

\author[P. Thieullen]{Phillipe Thieullen}
\address{Department of Mathematics \\
         Univ. Paris-Sud \\
         91405 Orsay, Cedex, France.}
\email{thieu@@topo.math.u-psud.fr}
\begin{document}

  \begin{abstract}
          Let $\s:\SS\hookleftarrow$ be a one-sided subshift
          of finite type. We show that for a generic $\a$-H\"older
          continuous function $A:\SS\to\re$, the supremum
          $$
          m(A)=\sup\{\textstyle\int A\,d\nu\,\vert\, \nu
          \text{ is a $\s$-invariant Borel probability }\}
          $$
          is achieved by a unique invariant probability.
          In the set $\cup_{\be>\a}C^\be(\SS,\re)$, with
          the $C^\a$-topology, generically the maximizing 
          measure is supported on a periodic orbit. This
          proves a version of a conjecture of R. Ma\~n\'e.
          We also show that these maximizing measures
          can be obtained as weak limits of equilibrium states.
          
          We apply these theorems to the class $\cF_\la(\a)$ of 
          $C^{1+\a}$ endomorphisms of the circle $f:S^1\to S^1$ which
          are coverings of degree 2, expanding $f'(x)>\la>1$, $\forall
          x\in S^1$ and orientation preserving. We prove that
          generically on $f\in\cF_\la(\a)$, the invariant probability
          which maximizes the Lyapunov exponent $\int \log f' \;d\nu$
          is unique, and that on $\cup_{\be>\a}\cF_\la(\be)$ (with
          the $C^{1+\a}$-topology) this (unique) maximizing measure 
          is suported on a periodic orbit.
 \end{abstract}

 \maketitle

\bigskip
\bigskip

 \noindent{\large\bf Introduction.}
 
 \bigskip
 \bigskip

 Let $\s:\SS\hookleftarrow$ be a one-sided topologically mixing
 subshift of finite type and $A:\SS\to\re$ a H\"older continuous 
 function. In this paper we are interested on $\s$-invariant 
 probability measures $\mu$ which maximize the integral 
 $\int A\, d\nu$ among all Borel $\s$-invariant probability 
 measures.

 Fix $0<\la<1$ and endow $\SS$ with the metric
 $d({\mathbf x},{\mathbf y})=\frac  1{\la^n}$, where
 ${\mathbf x}$,~${\mathbf y}\in\SS$,
 ${\mathbf x}=(x_0,x_1,\ldots)\in\SS$,
 ${\mathbf y}=(y_0,y_1,\ldots)$, and 
 $n=\min\{k\ge 0\,\vert\, x_k\ne y_k\}$. 
 Let $0<\a\le 1$, given an $\a$-H\"older function $A:\SS\to\re$, write
   $$
    \Hold_\a(A)
     =\sup_{0<d(x,y)\le 1}\frac{\abs{A(x)-A(y)}}{d(x,y)^\a},
     \qquad
   \norm{A}_0 = \sup_{x\in \SS}\abs{A(x)}  
   $$
   and define the $\a$-H\"older norm of $A$ by
   $$
   \norm{A}_\a = \Hold_\a(A)+\norm{A}_0.
   $$
   Denote by $\Can$ the set of $\a$-H\"older continuous
   functions $A:\SS\to\re$ endowed with the $\a$-H\"older norm
   $\norm{\;\;}_\a$. In view of applications, we shall restrict
   ourselves to H\"older functions with zero topological pressure.
   The results below hold also without this restriction. Denote
   by $\Ca$ the subset of functions $A\in\Can$ which have zero
   topological pressure. We shall prove
   
   \bigskip
   
   \noindent{\bf Theorem A.}
   
   {\it There is a generic set $\cG_1\subseteq\Ca$ such that
        if $A\in\cG_1$ then $A$ has a unique maximizing measure
        whose support is uniquely ergodic.   
   }     
   
   \bigskip

   The problem we consider here is in some sense analogous
   to the problems considered in the Aubry-Mather theory
   (see \cite{FM}, \cite{Ma2}) for Lagrangian flows.
   In particular our result is analogous to a recent result
   of Ma\~n\'e~\cite{Ma1} on Lagrangian flows, where he shows
   that generically (on the Lagrangian and on the homological
   position) there is a unique action minimizing measure.
   The main difference among these theories is that in the
   lagrangian setting the dynamics is defined by variational
   properties and hence minimizing properties imply
   invariance under the Euler-Lagrange flow. In our setting
   we have to impose somehow the invariance under the shift.
   The analogous to fix the homology in the Aubry-Mather
   theory in our setting is to consider side conditions,
   like $\int\psi_i\, d\nu=c_i$,
   $i=1,2,\ldots,k$, where $\psi_i\in\Can$ and $c_i\in\re$ are
   constants, in the maximization problem for $A$.       
   We obtain the same results in this case because
   by means of the Legendre transform this problem is equivalent
   to maximizing $A+\sum_{i=1}^k x_i\, \psi_i$ for certain fixed 
   values of $x_i\in\re$ (which depend on the $c_i$'s).
   
   In~\cite{Ma2} and~\cite{Ma3}, R. Ma\~n\'e conjectured that 
   generically the unique minimizing measure is supported on a
   periodic orbit. In our case we can prove this conjecture 
   in a subset of $\Ca$ of functions which have  slightly more
   regularity. 
     Let $\Cap$ be the closure in the $\a$-H\"older topology of
  $\cup_{\ga>\a}C^\ga_0(\SS,\re)$. Given a periodic point
  $p\in\Fix\,\s^n\subset\SS$, let $\nu_p$ be the $\s$-invariant
  probability supported on the positive orbit of $p$.
  
  \bigskip
  
  \noindent{\bf Theorem B.}
  
  {\it Let $\cG_2\subset\Cap$ be the set of $A\in\Cap$ such that
           there is a unique minimizing measure which is supported
           on a periodic orbit and it is locally constant. 
           Then $\cG_2$ is open and dense in $\Cap$. 
  }         
           
  \bigskip         
           
             Here locally constant means that             
             there is a periodic point $p\in\SS$ and
              a neighbourhood $\cU\ni A$ such that
             for all $B\in\cU$, the unique maximizing measure 
             for $B$ is $\nu_p$, where $\nu_p$ is the $\s$-invariant
             probability supported on the positive orbit of $p$.

  \bigskip
  
  The techniques used to prove this theorem involve the analogous to
  the {\it action potential} defined by Ma\~n\'e for Lagrangians
  in~\cite{Ma3}. Here we define this potential by
  $$
  S(x,y) := \textstyle\lim_{\e\to 0}\;\sup\big\{
  \sum\limits_{k=0}^n \bigl[A(\s^k z)-m_0\bigr]\,\big\vert\,
  n>0,\;\s^nz=y,\;d(z,x)<\e\,\big\},
  $$
  for $x,\;y\in\SS$, where $m_0=\inf\{\,\int A\,d\nu\,\vert\,
  \nu\text{ a $\s$-invariant Borel probability}\,\}$.
   In general, the function $S(x,y)$ is highly
  discontinuous (cf. proposition~\ref{P:discontinuous}), but
  if e.g. $x$ is in the support of a maximizing measure,
  then the function $y\mapsto S(x,y)$ is $\a$-H\"older continuous.
  Writing $V(y)=S(x,y)$ in this case, it is staright forward from
  the definition of $S$ that
  \begin{equation}\label{E:subsolution}
  V(\s y) \ge V(y) + A(y) - m_0
  \end{equation}          
  for all $y\in\SS$. In the Lagrangian case this $V$ corresponds
  to the existence of a subsolution of the Hamilton-Jacobi equation
  (cf. Fathi~\cite{Fa1}, \cite{Fa2}, also~\cite{CIPP}). Writing 
  $B(y)=A(y)+V(y)-V(\s y)$, then $B(y)$ is $\a$-H\"older and
  $\int B\,d\nu= \int A\,d\nu$ for any $\s$-invariant probability.
  Hence we can replace $A$ by $B$ in the maximization problem,
  with the advantage that $B\le m_0= \max_\nu \int B\, d\nu$.
  This implies that the inequality~\eqref{E:subsolution} in in fact
  an equality on the support of any minimizing measure. This, in turn,
  implies the
  
  \bigskip
  
  \noindent{\bf Coboundary Property.}
  
  {\it The function $A$ is cohomologous to a constant on the support
  of any maximizing measure by (the same) a H\"older continuous
  coboundary function, i.e. $A=V - V\circ\s+m_0$
  on $\supp(\mu)$ for any maximizing measure $\mu$.
   In particular any any measure supported 
  on a support of a maximizing measure is maximizing.
  }
  
  \bigskip
  
  In fact, the coboundary property can be extended to the set
  \linebreak
  $\fS =\{\, x\in \SS\,\vert\, W(x,x)=0\,\}$, which contains
  the support of all maximizing measures (cf.
  proposition~\ref{P:G.3} item [4]).

  It is possible to construct examples in  which there is
  a unique maximizing measure with positive entropy. In particular
  not supported on a periodic orbit. 
  
  There is an example in~\cite{G}
  of an invariant measure $\mu$ on the full 2-shift 
  $\Sigma^+_2=\{0,1\}^\na$ whose support is uniquely ergodic. 
  If $A:\Sigma^+_2\to\re$ is an $\a$-H\"older function which 
  attains its maximum exactly at $\supp(\mu)$, then $\mu$ 
  is the unique maximizing measure for $A$. By adding a 
  constant we can make $P(A)=0$.

  Another important ingredient in the proof of theorem B is
  the continuously varying support property, that we state now.
  Let $A\in\Ca$ and $\mu$ a maximizing measure for $A$.
  We say that a sequence of probability measures $\nu_n$
  {\it strongly converges\/} to a probability $\mu$
  if $\nu_n\to\mu$ weakly* and $\supp(\nu_n)\to\supp(\mu)$ in
  the Hausdorff metric. We say that the pair $(A,\mu)$ has the 
  {\it continuously varying support  property\/} if given a
  sequence $A_n\to A$ in $\Ca$,  an maximizing measures $\mu_n$ 
  for $A_n$, then $\mu_n$ strongly converges to $\mu$.
  
  \bigskip
  
  \noindent{\bf Theorem C.}
  
  {\it  
    There is a dense subset $\cD\subset\Cap$
    such that any $A\in\cD$ has a unique maximizing measure
    $\mu$ and the pair $(A,\mu)$ has the continuously varying
    support property. 
 }
 
 \bigskip

  Finally, we relate our maximization problem with the thermodynamic
  formalism. Let $A\in\Can$ and for each $t\in\re$ let $\hmu_t$ be the
  equilibrium state for $t\,A$. The following proposition appeared
  in a slightly different form in~\cite{L}:
  
  \bigskip

  \noindent{\bf Proposition D.}
  
  {\it
  Suppose that the maximizing measure $\mu_A$ for $A\in\Can$
  is unique and $A>0$. Then $\mu_A=\lim\limits_{t\to +\infty}\hmu_t$ in the
  weak* topology.
  }
  
  \medskip
  
  Recall the variational principle for the topological pressure
  $P(tA) = \max_\nu h(\nu)+t\int A\, d\nu$, where the maximum is
  along the $\s$-invariant probabilities, $P(tA)$ is the topological
  pressure and $h(\nu)$ is the metric entropy of $\nu$. 
  The result above shows that when $t\to +\infty$ in the variantional
  principle, one is putting less strength in the entropy of the
  measure and more stress in  the integral of the measure.
  However, the integral-maximizing measures do not seem to
  inherit properties from the approximating equilibrium states.

  \bigskip
  \bigskip
  
  \noindent{\bf Expanding maps of the circle.}
  
  \bigskip

  We can apply the results above to concrete situations
  using symbolic dinamics. An example that motivated us 
  is the case of invariant probabilities maximizing
  the Lyapunov exponent on a degree 2 expanding maps of the
  circle.

 Consider the class $\cF=\cF_\la(\a)$ of $C^{1+\a}$ endomorphisms of the
 circle $f:S^1\to S^1$ which are coverings of degree 2, expanding
 $f'(x)>\la>1$, $\forall x\in S^1$ and orientation preserving.
 For a $C^{1+\a}$ endomorphism $f\in\cF$, denote its $C^{1+\a}$
 norm by
 $$
 \norm{f}_{1+\a} = \norm{f}_{C^1}+\Hold_\a(f').
 $$
 
 We say that an $f$-invariant Borel probability is a {\it Lyapunov
 maximizing measure\/} or simply a {\it maximizing measure\/}
 if it maximizes the integral
 \begin{equation}\label{E:Lya}
  \int \log f'\;d\nu 
 \end{equation}
 among all $f$-invariant probabilities. The Lyapunov exponent of an
 invarinat measure $\nu$ is given by the integral~\eqref{E:Lya}. It
 expresses the mean value of the rate of expansiveness of points in
 the support of $\nu$. We are looking for measure with maximal
 sensitivity dependence on initial conditions.
 
 \bigskip
 
 \noindent{\bf Theorem A1}
  
 {\it Generically on the $C^{1+\a}$-topology for maps $f\in\cF$, there
 exists a unique $f$-invariant Lyapunov maximizing measure $\mu_f$. 
 Moreover, the support of $\mu_f$ is uniquely ergodic.  }
 
 \bigskip

 Let $\cF(\a+)$ be the closure of $\cup_{\be>\a} \cF_\la(\be)$ 
 in $\cF$ (with the $C^{1+\a}$-topology).
 
 \bigskip  
 
 \newpage
 
 \noindent{\bf Theorem B1}
 
 {\it There is a generic set $\cG_2\subset\cF(\a+)$ such that for
  $f\in\cG_2$ there is a unique $f$-invariant Lyapunov maximizing
  measure and it is supported on a periodic orbit.
  }
  
 \bigskip
 \bigskip
 
 In section~\S\,\ref{UNIQUENESS} we prove theorem~A. On
 section~\S\,\ref{SHADOWING} we show some preliminary shadowing lemmas.
 On section~\S\,\ref{POTENTIAL} we define and prove the properties 
 of the action potential and state the coboundary property. 
 On section~\S\,\ref{CVS} we prove theorem~C. On
 section~\S\,\ref{PERIODIC} we prove theorem~A. On
 section~\S\,\ref{THERMO} we prove theorem~D. On section~\ref{CIRCLE}
 we prove theorems~A1 and B1 and give an equivalence between
 $C^{1+\a}$ expanding dynamics on the circle and $\a$-H\"older maps
 on the shift.
 
 \bigskip
 \bigskip

  \refstepcounter{section}\label{UNIQUENESS}
  \noindent{\large\bf \thesection. Generic uniqueness of the
                             maximizing measure.}

  \bigskip                           
  \bigskip

   In this section we prove theorem A.
   We start with nome notation. Denote by $K(\s)$
   the set of $\s$-invariant Borel probabilities on $\SS$, endowed
   with the weak* topology. Given $A\in\Ca$, set
   \begin{align*}
   m(A) &= \max\bigl\{\,\textstyle \int A\, d\nu\,\big\vert\,
            \nu\in K(\s)\,\bigr\} \\
   \cM(A) &= \bigl\{\, \mu\in K(A)\,\big\vert\,
            \textstyle\int A\,d\mu = m(A)\,\bigr\}
   \end{align*}                 
   A measure $\mu$ in $\cM(A)$ is called a {\it maximizing measure\/}.
   
   The arguments in this section rely on Banach space techniques.
   Since $\Ca$ is a Banach manifold, we need to do some conversions:
   
   \bigskip
   
   \noindent{\bf Proof of theorem A:}
   
   The topological pressure $P:\Can\to\re$ is real analytic
   (cf.~\cite{PP}). It is also a submersion because $P(A+r)=P(A)+r$
   for any $r\in\re$. Hence the set $\Ca=P^{-1}\{0\}$ is a Banach
   manifold. 
   
   Given $A_0\in\Ca$, the derivative of the pressure at $A_0$ is
   given by $DP(A_0)\cdot B = \int B\, d\hmu$, where $\hmu$ is
   the equilibrium state for $A_0$ and $B\in\Can$ (see~\cite[corollary
   1.4 and 1.7]{Ma1} or~\cite{PP} ).
   Hence the tangent space at $A_0$ to $\Ca$ is the set of functions
   $A_0+B$ where $B\in\Can$ and  $\int B\, d\hmu =0$. This space does
   not intersect $0$ because $\int A_0\, d\hmu = -\text{entropy of
   }\hmu$.
   
  Locally, near $A_0$, there is a homeomorphism $\cY$ associating
  $A\in\Ca$ to $A_0+B$ $B\in T_{A_0}\Ca$. This homeomorphism is of
  the form $\cY(A)=c_A\, A = A_0 + B$, where $c_A$ is a positive
  constant. Therefore the maximizing measures for $A$ or $A_0+B=c_A\,
  A$ are the same. Thus to show the generic properties of maximizing
  measures for $A$ or $A_0+B=c_A\, A$ is the same problem.
 
  So, we have to show that generically on functions $B$ close
  to zero and such that 
  \linebreak
  $\int B\, d\hmu=0$, the function
  $A_0+B$ has a unique maximizing measure.
  This is proven on theorem~\ref{uniqueness} below.
  
  \qed

  \bigskip
  
  Fix a Borel probability measure $\hmu$ on $\SS$ and $A_0\in\Ca$.
  Define
  $$
  \cH = \{ \, A_0+B\,\vert\, \textstyle\int B\, d\hmu = 0 \,\}
  $$
  
  \bigskip
  \bigskip

 \begin{Theorem}\label{uniqueness}
 There exists a residual set $\cG_1\subset \cH$
 such that for all $A\in\cG_1$, the set $\cM(A)=\{\mu\}$ 
 has a unique element. Moreover the support of $\mu$ is uniquely
 ergodic.
 \end{Theorem} 
 
 \bigskip
 
 \noindent{\bf Proof:}
 
 The proof will require two lemmas. The idea is to show that for
 any $\e>0$, the open set
 $$
 \cO_\e=\{\, A\in\cH\,\vert\,\diam\cM(A)\le\e\,\}
 $$
 is dense. Considering $\e$ of the form $\e=\frac 1n$, $n\in\na$, 
 we obtain from Baire's Theorem that there is a residual set
 whith a unique maximizing measure. Item~[4] of
 proposition~\ref{P:G.3}
 implies that any invariant measure on the support of a maximizing 
 measure is maximizing. Hence if there is a unique maximizing measure,
 its support must be uniquely ergodic.

 Consider $K_0$ a subset of $K(\s )$ (the set of invariant measures for
 $\s $) and define 
 \begin{gather*}
 m_0(A) := \sup \bigl\{\, \textstyle\int A\, d\nu\,\vert\,\nu\in
 K_0\,\bigr\}, \\
 \cM_0(A):=\bigl\{\,\mu\in K_0\,\vert\,\textstyle\int A\,
 d\mu=m_0(A)\,\bigr\}.
 \end{gather*}
 
 We say that $u$ is an extremal point of the convex set $C$, if $u$ is
 not a mid point of a segment where the endpoints are in $C$. A point
 $p$ in the convex set $C$ is said a strictly extremal point for $C$,
 if there exists a linear map $h$ on the set $E$ such that the
 supremum of $h$ restricted to $C$ is attained at $p$ and only at $p$.
 
 A classical result in convex analysis (see Strazewicz's Theorem in
 ~\cite{R}) states that any extremal $u$ can be approximated by a
 strictly extremal $p$.
 
 \bigskip
 
 \begin{Lemma}\label{l.6}\quad
 
If $\mu_0$ is an extremal point of a compact set $K_0$, the for all
$\e>0$, there exists $w\in \cH$ such that $\diam\cM_0(w)<2\,\e$ and
$d(\mu_0,\cM_0(w))<\e$.
\end{Lemma}

\medskip

\noindent{\bf Proof:} 

Consider a sequence $w_n$,~$n\in\na$, of functions in $\cH$ that
define a metric $\td$ on $K_0$ by
$$
\td(\nu,\mu)=\zum_{j=1}^\infty 
\displaystyle \frac 1{2^j}\,
\Big|\int w_j\, d\mu - \int w_j\,d\nu\,\Big|,
$$ 
compatible with the weak convergence on the compact space of
probabilities on $\SS$. For each $n\in\na$, define $P_n:K_0\to\re^n$
by 
$$
P(\mu)=\big(\textstyle\int w_1\,d\mu,\ldots,\textstyle\int w_n\,
d\mu\big).
$$
From the definition of $\td$ and the compactness of $K_0$, it is easy
to see that for all $\e>0$, there exist $\de>0$ and $n>0$ such that,
if
$S\subset\re^n$ and $\diam \,S<\de$, then 
\begin{equation}\label{E:a.1}
\diam(P_n^{-1})(S)<\e.
\end{equation}

Note that $u=P_n(\mu_0)$ is an extremal point of $P_n(K_0)=C$. From
Strazewicz's Theorem applied to $C$, let $p$ be a strictly extremal
point such that $d(p,P_n(\mu_0))<\de$.  Then by~\eqref{E:a.1}, we have
that
\begin{equation}\label{E:a.2}
\diam(P^{-1}(p),\mu_0)<\e.
\end{equation}
Consider $w=\sum_{i=1}^n\la_i\,w_i$, where the $\la_i$ are such that
$h(x_1,x_2,\ldots,x_n)=\sum_{i=1}^n\la_i\,x_i$. It is easy to see that
$\cM_0(w)=P^{-1}(p)$. Therefore, by~\eqref{E:a.2},
$\td(\mu_0,\cM_0(w))<\e$ and $\diam(\cM_0(w))<2\,\e$.
This shows the lemma.

\qed

\bigskip

\begin{Lemma}\label{l.7}\quad

Suppose that $A\in\cH$ and consider and extremal point $\mu_0$ of
$\cM(A)$. Then for any neighbourhood $U$ of $\mu_0$ in $K(\s )$, and
every $\e>0$, there exists $A_1\in U$ such that $d(\mu_0,\cM(A_1))<\e$
and $\diam\cM(A_1)<2\,\e$.
\end{Lemma}

\medskip

\noindent{\bf Proof:}

 Applying lemma~\ref{l.6} to $K_0=\cM(A)$, for any $\e>0$ there exists
 $w$ such that $\diam\cM_0(w)<\e$ and $d(\mu_0,\cM_0(w)<\e$. Let
 \begin{align*}
 m  &=m_0(w)=\sup\bigl\{\,\textstyle\int
       w\,d\mu\,\big|\,\mu\in\cM(A)\,\bigr\},  \\
 m_0&= m(A) =\sup\bigl\{\,\textstyle\int A\,d\mu \,\big\vert\,
      \mu\in K(\s )\,\bigr\}.
 \end{align*}
 
 Denote by $f_0$ and $f_1$ the functions defined on $\mu\in K(\s )$ by
 \begin{align*}
 f_0(\mu) &= \int A\, d\mu - m_0, \\
 f_1(\mu) &= \int w\, d\mu - m .
 \end{align*}
 Then 
 \begin{align}
 f_1(\mu) &=0, \qquad\text{ for all } \mu\in\cM_0(\mu),\label{e.5} \\
 f_0(\mu) &=0, \qquad\text{ for all } \mu\in\cM_0(\mu),\label{e.6}
 \end{align}
 (because $\cM_0(w)\subset K_0=\cM(A)$).
 
 Reciprocally, observe that if $\mu\in\cM(A)=K_0$ and if 
 $f_1(\mu)=0$, then 
 \begin{equation}\label{e.7}
 \mu\in\cM_0(w).
 \end{equation} 
 For $\mu\in K(\s )$, if
 $f_0(\mu)=0$, then 
 \begin{equation}\label{e.8}
 \mu\in K_0=\cM(A).
 \end{equation}
 
 Observe that by the definition of $m$ and $m_0$,
 \begin{equation}\label{e.9}
 f_1(\mu)\le 0 \quad\text{ for all }\mu\in K_0
 \quad\text{ and }\quad 
 f_0(\mu)\le 0 \quad\text{ for all }\mu\in K(\s ).
 \end{equation}
 
 Now define$f_\la=f_0+\la\, f_1$ for all $\la>0$. 
 Let 
 \begin{gather*}
 m(\la) = \max_{\nu\in K(\s )}f_\la(\nu),\\
 \cM_\la =\{\,\mu\in K(\s )\,\vert\,f_\la(\mu)=m(\la)\,\}.
 \end{gather*}
 Observe that 
 \begin{equation}\label{e.10}
 m(\la)\ge 0,
 \end{equation}
 because by~\eqref{e.5} and~\eqref{e.6},
 if $\mu\in\cM_0(\mu)$, then $f_\la(\mu)=0$.
 
 \bigskip
 
 \noindent{\bf Claim:} $\lim_{\la\to 0}\;\diam(\mu_0,\cM(\la))<\e$.
 
 \medskip
 
 If this claim is true,  taking $A_1=f_\la$ for $\la$ small
 the lemma is proved.
 
 Suppose that the claim is false. The there exist a sequence
 $\la_n\to 0$ and $\mu_n\in\cM(\la_n)$ such that
 \begin{equation}\label{e.11}
 \inf_{n\in\na}\td(\mu_n,\mu_0)\ge \e.
 \end{equation}
 Consider a limit $\ov{\mu}$ of a subsequence of $\mu_n$.
 Then by~\eqref{e.11} 
 \begin{equation}\label{e.12}
 \td(\ov{\mu},\mu_0)\ge\e.
 \end{equation}
 If we prove that $\ov{\mu}\in\cM_0(w)$, then from~\eqref{e.12},
 we obtain a contradiction.  Note that
  \begin{equation}\label{e.13}
 f_1(\mu_n)\ge 0,
 \end{equation}
 because from~\eqref{e.10} and~\eqref{e.9}
 $$
 0\le m(\la_n)=f_0(\mu_n)+\la\, f_1(\mu_n)\le \la_n\, f_1(\mu_n).
 $$
 Note also that 
 $$
 \lim_{n\to\infty}m(\la_n)=\lim_{n\to\infty}\max_{\mu\in
 K(\s )}f_{\la_n}(\mu)=\max_{\mu\in K(\s )}f_0(\mu)=0.
 $$
 Since $f_{\la_n}=f_0(\mu_n)+\la_n \, f_1(\mu_n)=m(\la_n)$,
 by continuity,
 $$
 f_0(\ov{\mu})=f_0\bigl(\lim_{n\to\infty}\mu_n\bigr)
             =\lim_{n\to\infty}\bigl(m(\la_n)-\la_n\, f_1(\mu_n)\bigr)
             =0.
 $$
 Therefore by~\eqref{e.8}, $\ov{\mu}\in\cM(A)$. Now, by~\eqref{e.13},
 $f_1(\mu_n)\ge 0$ and then $f_1(\ov{\mu})=\lim_{n\to\infty}\ge0$.
 Since $\ov{\mu}\in K_0=\cM(A)$, then by~\eqref{e.9},
 $f_1(\ov{\mu})\le 0$. Therefore $f_1(\ov{\mu})=0$.
 Finally, from~\eqref{e.7} $\ov{\mu}\in\cM_0(A)$. This contradicts
 $\td(\ov{\mu},\mu_0)\le \e$.
 \qed

 \bigskip

  \bigskip
  \bigskip
  
 \refstepcounter{section}\label{SHADOWING}
 \noindent{\large \bf \thesection. Shadowing Lemmas.}

  \bigskip
  \bigskip

 Let $\s:\SS\hookleftarrow$ be a positive subshift of finite type.
 For ${\mathbf x}$,~${\mathbf y}\in\SS$, 
 ${\mathbf x}=(x_0,x_1,\ldots)\in\SS$,
 ${\mathbf y}=(y_1,y_2,\ldots)$, 
 write $d(\mathbf{x},\mathbf{y})=\frac 1{\la^n}$, 
 where $n=\min\{k\ge 0\,\vert\, x_k\ne y_k\}$. Let 
 $\e_0=\frac 1\la>0$, so that
 \begin{itemize}
 \item[{\rm [1]}] If $y=(y_0,y_1,\ldots)\in\SS$ and 
            $(x_0,y_0,y_1,\ldots)\in\SS$ 
            then the local inverse
            $\psi_{x_0}(z)=(x_0,z_0,z_1,\ldots)$ is defined on all 
            $\{\, z\,\vert\, d(z,y)<\e_0\,\}$.
 \item[{\rm [2]}] If $x$,~$y\in\SS$ and $d(x,y)<\e$, then
            $d(\sigma x,\sigma y)=\la\, d(x,y)$.
 \end{itemize}
 
 \bigskip
 
 For $x\in\SS$ and $r>0$ write 
 $B(x,r)=\{\, y\in\SS\,\vert\,d(x,y)<r\,\}$.

 \medskip
 
 We say that a sequence $\{x_0,\ldots,x_n\}$ is a 
 {\it $\de$-pseudo-orbit with $M$ jumps}, 
 if $d(\s x_i,\s x_{i+1})\le\de$ for all $0\le i\le N-1$ and 
 $\#\{\,0\le i\le N-1\,|\,x_{i+1}\ne\s(x_i)\,\}=M$.
 We say that a $\de$-pseudo-orbit $\{\,x_0,\ldots,x_N\,\}$ is
 {\it $\e$-shadowed by $p\in\SS$}, 
 if $d(\s^kp,x_k)<\e$ for all $0\le k\le N$.
 
 \bigskip

 \begin{Lemma}\label{L:G.2}
 Let $\e_1:=(1-\la^{-1})\,\e_0$. For all $A:\SS\to\re$ $\a$-H\"older
 continuous, there exists $K_1=K(A,\la)>0$ such that if $0<\de<\e_1$ and
 $\{\,x_0,\ldots,x_N\,\}$ is a $\de$-pseudo orbit with $M$ jumps,
 then there exists $p\in\SS$ that $\big(\frac \de{1-\la^{-1}}\big)$-shadows
 $\{x_i\}_{i=1}^N$ and for all $0\le i\le j\le N$
 \[
 \bigg\vert{ \textstyle
        \sum\limits_{k=i}^{j} A(\s^kp)-\sum\limits_{k=i}^j A(x_k)}
 \bigg\vert
 \le M\, K_1\, \de^\a.       
 \]
 
 Moreover, 
 \begin{itemize}
 \item[{\rm [1]}] The point $p$ can be taken such that $\s^N(p)=x_N$.
 \item[{\rm [2]}] If the pseudo-orbit is periodic (i.e. $x_N=x_0$),
            then the point $p$ can be taken $N$-periodic:
            $\s^N(p)=p$.
 \end{itemize}           
 \end{Lemma}       
 
 \medskip

 \medskip

 \noindent{\bf Proof:}
 
 For $1\le n\le N$ let $\vr_n:B(x_n,\e_0)\to\SS$ be the branch of the
 inverse of $\s$ such that $\vr_n(\s x_{n-1}) =x_{n-1}$. Then
 $\psi_N:=\vr_1 {}_\circ\vr_2 {}_\circ\ldots {}_\circ \vr_N$ is a
 contraction with Lipschitz constant $\la^{-N}<1$. Moreover,
 $\vr_n(B(x_n,r))\subseteq\vr_n(B(\s x_{n-1}, r+\de))\subseteq
 B(x_{n-1},r)$ for $r=\frac{\de}{\la-1}$, $r+\de<\e_0$. [This gives
 $\de<(1-\la^{-1})\,\e_0 =:\e_1$.] In particular
 $\psi_N(B(x_N,r))\subseteq B(x_0,r)$.
 \begin{itemize}
 \item[{\rm [1]}] Let $p=\psi_N(x_N)\in B(x_0,r)$.
 \item[{\rm [2]}] Let $p\in B(x_0,r)=B(x_n,r)$ be the fixed point of
            $\psi_N$.
 \end{itemize}
 Then $d(\s^k p,x_k)\le r = \frac{\de}{\la -1}$.
 
 Let $0<a_1<a_2<a_3<\cdots<a_M\le N$ be the indices such that
 $\s(x_{a_i})\ne x_{a_i+1}$. Let $a_0=0$, $a_{M+1}=N$ and
 $b_i=a_{i+1}-a_i$, $0\le i \le M$. Then, for $0\le j<b_i$, we have
 that
 \begin{align*}
 d(\s^{a_i+j} p, x_{a_i+j}) 
   &\le \la^{j-b_i}\; d(\s^{a_{i+1}} p, \s^{b_i} x_{a_i} ) \\
   &\le \la^{j-b_i}\;\bigl[d(\s^{a_{i+1}}p,x_{a_{i+1}})
        + d(\s(x_{a_{i+1}-1}),x_{a_{i+1}})\bigr]   \\
   &\le \la^{j-b_i}\; \left[\tfrac{\de}{\la-1} + \de\right]
   = \la^{j-b_i}\; \tfrac 1{1-\la^{-1}}\; \de .
 \end{align*}
 \begin{align*}
 \Big\vert\zum_{k=i}^j A(\s^k p) - \zum_{k=i}^j A(x_k)\Big\vert
 &\le \zum_{k=0}^{N-1}\abs{A(\s^k p)-A(x_k)} \\
 &\le \zum_{k=0}^M \zum_{j=0}^{b_k-1}
      \Hold_\a(A)\; \la^{-j\a}\;\tfrac{1}{(1-\la^{-1})^\a}\; \de^a \\
 &\le (M+1)\; \Hold_\a(A)\; \tfrac 1{1-\la^{-\a}}\;
      \tfrac 1{(1-\la^{-1})^\a}\; \de^\a.
 \end{align*}            
 \qed

 \bigskip 
 \bigskip

 \refstepcounter{section}\label{POTENTIAL}
 \noindent{\large\bf \thesection. The Action Potential.}
 
 \bigskip
 \bigskip
 
 Given $x$,$y\in\SS$ and $\de>0$, define
 $$
 S_\de(x,y):=\sup\bigg\{\zum_{k=0}^{n-1}\b[A(\s^k(z)-m_0\b]\,\bigg|\,
 \s^nz=y,\;d(z,x)<\de\,\bigg\},
 $$
 where $m_0=m(A)$.
 Since $\SS$ is topologically transitive, the backward orbit of any
 point $y\in\SS$ is dense in $\s$. Hence the set in the definition
 above is non-empty and thus $S_\de(x,y)>-\infty$ for any $\de>0$.
 We will show below that $\sup_{x,y\in\SS,\de>0}S_\de(x,y)<+\infty$.
 Since the function $\de\mapsto S_\de(x,y)$ is increasing, we can define
 $$
 S(x,y)=\lim_{\de\to 0^+}S_\de(x,y).
 $$
 We get that $-\infty\le S(x,y)\le Q$. In fact the value
 $S(x,y)=-\infty$ is possible and in general the function
 $S(x,y)$ is highly discontinuous. We quote the properties of $S(x,y)$
 in the following proposition:
 
 \bigskip
 
 \newpage
 
 \begin{Proposition}\label{P:G.3}\quad
 
 \begin{itemize}
 \item[{\rm [1]}] There is $Q>0$ such that $S_\de(x,y)<Q$ for all
            $x$~$y\in\SS$ and all $\de>0$.
 \item[{\rm [2]}] For all $x\in\SS$, $S(x,x)\le 0$.
 \item[{\rm [3]}] For all $x,\;y,\;z\in\SS$, $S(x,y)+S(y,z)\le S(x,z)$.
 \item[{\rm [4]}] Let
            $$
            \fS :=\{\,x\in\SS\,\vert\,S(x,x)=0\,\}.
            $$
            Then $\fS$ is closed and foward invariant.
            A measure $\mu$ is maximizing if and only if
            $\supp(\mu)\subseteq\fS$. In particular
            $\fS\ne\0$.
 \item[{\rm [5]}] If $x\in\fS$ then the function $W:\SS\to\re$,
            $W(y)=S(x,y)$ is finite and $\a$-H\"older continuous
            with $\Hold_\a(W)\le C(\la)\, \Hold_\a(A)$.
            Moreover, $W(y)-W(x)\ge S(x,y)$ for all $x,\;y\in\SS$.
 \end{itemize}           
 \end{Proposition}
  
 \bigskip
 
% \newpage
 \begin{Corollary}\label{Cor}\quad
 \begin{itemize}
  \item[{\rm [1]}] The $\a$-H\"older continuous function 
            $B(x):= A(x)-m_0+W(x)-W(\s x)$
            satisfies $B\ge 0$, $\int B\, d\nu = \int A\, d\nu$ for 
            any invariant measure and $\int B\, d\mu =0$ for any
            maximizing measure.          
  \item[{\rm [2]}] If $\mu$ is a maximizing measure, then any 
            invariant measure $\nu$ with
            $\supp(\nu)\subseteq\supp(\mu)\; (\subset\fS)$
            is maximizing. In particular if $A$ has a unique
            maximizing measure, then the set $\fS$ 
            (and hence also $\supp(\mu)$) is uniquely ergodic.
 \end{itemize}
 \end{Corollary}
 
 \medskip
 
 \noindent{\bf Proof:}  
 Item~[1] follows from~\ref{P:G.3}~[5]
 because 
 $$W(\s x)- W(x)\ge S(\s x,x)\ge A(x).$$
  \qed
 
 \bigskip

 A subset $K\subseteq\SS$ is said {\it $\e$-separated\/} if
 $d(x,y)>\e$ for all $x,\;y\in K$ with $x\ne y$.
 Given a periodic point $p\in\Fix(\s^n)$, let $\nu_p$ be the
 probability measure defined by
 $$
 \int f\; d\nu_p = \frac 1n \zum_{k=0}^{n-1} f(\s^kx),
 $$
 for any continuous function $f:\SS\to\re$. 
%  Given a continuous
%  function $f:\SS\to\re$ define
%  $$
%  \lV f\rV_0 := \sup_{x\in\SS} \lv f(x)\rv.
%  $$
%  
 \bigskip
 
 \noindent{\bf Proof of proposition~\ref{P:G.3}:}
 
 [1] \quad Let $\e=\e_2$ from lemma~\ref{L:G.2}. Let
     $$
     M(\e):=\max\{\,B\subseteq\SS\,\vert\, B \text{ is
     $\e$-separated}\,\}.
     $$ 
     Let $N>M(\e)$ and $x\in\SS$. Let
     \begin{align*}
     k_0 &= \max\{\,0\le k\le N\,\vert\, \{x,\s x,\ldots,\s^kx\}
         \text{ is $\e$-separated}\,\},
         \\
     k_1 &= \max\{\,0\le k\le N\,\vert\, \{\s^kx,\s^{k+1}x,\ldots,\s^Nx\}
         \text{ is $\e$-separated}\,\}.
     \end{align*}    
     Then $k_0\le M(\e)$ and $N-k_1\le M(\e)$. The set
     $\{\s^jx\,|\,0\le j\le k_0,\;k_1\le j\le N\,\}$ is not
     $\e$-separated. Hence there are $0\le i\le k_0$, $k_1\le j\le N$
     such that $d(\s^ix,\s^jx)<\e$. By lemma~\ref{L:G.2},
     \begin{align*}
     \sum_{k=i}^j\b[A(\s^kx)-m_0\b]
       &\le \sup_{p\in\Fix\s^n} n{\displaystyle\int\b[A-m_0\b]\;d\nu_p
       + K\,\e^\a \le K\,\e^a.}  \\
     \intertext{and}
     \zum_{k=1}^{N-1}\b[A(\s^kx)-m_0\b]
       &\le K\,\e^\a+2\,M(\e)\,\norm{A-m_0}_0,
     \end{align*}
    for all $x\in\SS$ and {\it all} $N>0$.     
    Thus     
    $$
    S_\de(x,y)\le K\,\e^\a+2\, M(\e)\,\norm{A-m_0}_0
    $$
    for all $\de>0$, $x,\;y\in\SS$. This implies item[1].
    
 \bigskip   
    
 [2] \quad If $0<\de<\e_2$, $d(x,y)<\de$ and $\s^ny=x$, then by
    lemma~\ref{L:G.2},
    $$
    \lv\zum_{k=0}^{n-1}\b[A(\s^ky)-m_0\b]
       -\zum_{k=0}^{n-1}\b[A(\s^kp)-m_0\b]\rv
       \le K\, (2\de)^\a
    $$
    for some periodic point $p\in\Fix\;\s^n$. Since
    $\int[A-m_0]\;d\nu_p\le 0$, then
    $$
    \sum_{k=0}^{n-1} \b[A(\s^ky)-m_0\b]
    \le n\displaystyle \int\b[A-m_0\b]\;d\nu_p 
    + K\, (2\de)^\a \le K\, (2\de)^a.
    $$
    Hence $S_\e(x,y)\le K\, (2\de)^a = K\, (2(1-\la)\e)^a$
    for $e=\frac{\de}{1-\la}$. Letting $\e\to 0$, we obtain that
    $S(x,x)\le 0$.       
    
 \bigskip
 
 [3] \quad  Given $\de>0$ let $a,\;b\in\SS$ be such that $d(x,a)<\de$,
     $\s^na=y$; $d(y,b)<\de$, $\s^mb=z$ for some $n,\;m>0$ and
     \begin{align}
     \zum_{k=0}^{n-1}\b[A(\s^ka)-m_0\b] &\ge S_\de(x,y)-\de\label{E:3.1}
     \\
     \zum_{k=0}^{m-1}\b[A(\s^kb)-m_0\b] &\ge S_\de(y,z)-\de\label{E:3.2}
     \end{align}
     Then $\{a,\s(a),\ldots,\s^{n-1}a,b,\ldots,\s^mb=z\}$ is a
     $2\de$-pseudo-orbit with 1 jump. By lemma~\ref{L:G.2}, there is
     $p\in\SS$ which $\b[\frac{2\de}{1-\la}\b]$-shadows the pseudo-orbit,
     $\s^{n+m}p=\s^mb=z$ and
     $$
     \zum_{k=0}^{n+m-1}\b[A(\s^kp)-m_0\b]
     - \bigg[ \zum_{k=0}^{n-1}\b[A(\s^ka)-m_0\b]
              +\zum_{k=0}^{m-1}\b[A(\s^kb)-m_0\b]\bigg]
     \ge -K\, (2\de)^\a.
     $$
     Since $d(x,p)\le d(x,a)+d(a,p)\le \de\,\big[\frac
     2{1-\la}+1\big]=:\e(\de)$, and $\s^{n+m}p=z$, then,
     using~\eqref{E:3.1} and~\eqref{E:3.2}, we have that
     $$
     S_{\e(\de)}(x,z)
     \ge[S_\de(x,y)-\de]+[S_\de(y,z)-\de]-2^\a\, K\,\de^\a.
     $$
     Letting $\de\to 0$, then $\e(\de)\to 0$ and
     $$
     S(x,z)\ge S(x,y)+S(y,z).
     $$
                       
   \bigskip
   
   In order to prove item~[5] we need the following
   
   \bigskip
   
   \begin{Lemma}\label{L:3.4}
   If $S(x,x)=0$, then for all $\e>0$ and $M>0$ there exists
   $w\in\SS$ and $n>M$ such that $d(w,x)<\e$, $\s^nw=y$ and
   \begin{equation}\label{E:L.3.4}
   \zum_{k=0}^{n-1}\b[A(\s^kw)-m_0\b]\ge S(x,y)-\e.
   \end{equation}
   \end{Lemma}
   
   \medskip
   
   \noindent{\bf Proof:}
   
   Let $\de>0$ be such that
   \begin{align*}
   M\,\de+\de+M\,K\,\de^\a &<\e,   \\
   \tfrac{\de}{1-\la}+\de &<\e. 
   \end{align*}  
   Let $a\in\SS$ and $n>0$ be such that $d(a,x)<\de$, $\s^na=y$ and 
   $$
   \zum_{k=1}^{n-1}\b[A(\s^ka)-m_0\b]
   \ge S_\de(x,y)-\de \ge S(x,y)-\de.
   $$
   Since $S(x,x)=0$ then there is $b\in\SS$ and $n>0$ such that
   $d(b,x)<\de$, $\s^mb=x$ and 
   $$
   \zum_{k=1}^{m-1}\b[A(\s^kb)-m_0\b]\ge S_\de(x,x)-\de
   \ge S(x,x)-\de\ge -\de.
   $$
   The ordered set $\{b,\ldots,\s^{m-1}b\},\overset{\text{$M$
   times}}{\ldots\ldots\ldots}\{b,\ldots,\s^{m-1}b\},
   \{a,\ldots,\s^na\}$ is a $2\de$-pseudo-orbit with $M$ jumps.
   By lemma~\ref{L:G.2} there is $w\in\SS$ such that
   $d(b,w)<\frac{\de}{1-\la}$, $\s^{mM+n}w=\s^na=y$ and 
   \begin{align*}
   \bigg|\zum_{k=0}^{mM+n-1} &\b[A(\s^kw)-m_0\b]     \\
   &-\bigg\{M\zum_{k=0}^{m-1}\b[A(\s^kb)-m_0\b]
   +\zum_{k=0}^{n-1}\b[A(\s^ka)-m_0\b]\bigg\}\bigg|
   \le M\,K\,\de^\a.
   \end{align*}
   Then
   \begin{align*}
   \zum_{k=0}^{mM+n-1}\b[A(s^kw)-m_0\b]
     &\ge -M\,\de + S(x,y) -\de - M\,K\,\de^\a \\
     &\ge S(x,y)-\e.
   \end{align*}
   Moreover $d(w,x)\le d(w,b)+d(b,x)\le\frac{\de}{1-\la}+\de<\e$,
   $s^{mM+n}w=y$ and $mM+n>M$.
   
   \qed  
 
   \bigskip
   
   [5] \quad Now we prove item~[5]. Let $z,\;y\in\SS$ and $d(y,z)=d$
   small. Given $\e>0$ let $M=M(\e)>0$  be such that $\la^M(\e+d)<\e$.
   Let $w\in\SS^+$ and $n>M(\e)$ be as in lemma~\ref{L:3.4}. Since
   $d(\s^nw,z)\le d(\s^nw,y)+d(y,z)\le\e+d$, then the ordered set
   $\{w,\s w,\ldots,\s^{n-1}w,z\}$ is an $(\e+d)$-pseudo-orbit with 1
   jump. By lemma~\ref{L:G.2} there exists $p\in\SS$ such that
   $s^np=z$, $d(w,p)<\la^n (d+\e)$ and
   \begin{equation}\label{E:L.3.4.2}
   \bigg|\zum_{k=0}^{n-1}\b[A(\s^kp)-m_0\b]
        -\zum_{k=0}^{n-1}\b[A(\s^kw)-m_0\b]\bigg|
   \le K\, (d+\e)^\a.
   \end{equation}
   Since $n>M(\e)$ we have that $d(p,x)\le
   d(p,w)+d(w,x)\le\e+\la^n(d+\e)<2\e$. Then, using~\eqref{E:L.3.4.2}
   and~\eqref{E:L.3.4}, we have that
   $$
   S_{2\e}(x,z)\ge\zum_{k=0}^{n-1}\b[A(\s^kp)-m_0\b]
   \ge S(x,y)-\e-K\, (d+\e)^\a,
   $$
   where $K=C(\la)\,\Hold_\a(A)$.
   Letting $\e\to 0$ we get that 
   $$
   S(x,z)\ge S(x,y)-K\, d^\a.
   $$
   Interchanging the roles of $y$ and $z$ we obtain that
   $$
   \lv W(y) - W(z)\rv=\lv S(x,y)-S(x,z)\rv\le K\, d^\a.
   $$
   Now, by the triangle inequality, we have that
   $$
   W(z)-W(y)=S(x,y)-S(x,y)\ge S(y,z).
   $$
   
   \bigskip

   [4]  \quad We now prove item~[4]. We first prove that if $\mu$
   is an invariant measure with $\supp(\mu)\subseteq\fS$ then it is
   maximizing. Fix $x\in\fS$ and define $W(y)=S(x,y)$ and
   $B(y)=A(y)-m_0+W(y)-W(\s y)$. By item~[4] we have that
   $W(\s y)-W(y)\ge S(y,\s y)\ge A(y)$. Hence $B(y)\le 0$ for all
   $y\in\SS$ and $\int B\, d\mu = \int (A-m_0)\, d\mu$.
   
   To see that $\mu$ is maximizing, it is enough to show that
   $B\equiv 0$ on $\fS$. Let $y\in\fS$. Then $S(y,y)=0$ and for 
   any $\de>0$ there exists $z=z(\de)\in\fS$ and $n>0$ such that
   $d(z,y)<\de$, $\s^nz=y$ and
   $$
   \zum_{k=0}^{n-1}\b[A(\s^kz)-m_0\b] > S_\de(y,y)-\de.
   $$
   Then
   $$
   S_{\la\de}(\s y,y)\ge \zum_{k=1}^{n-1}\b[A(\s^k z)-m_0\b]
   > S_\de(y,y)-\de A(z)+m_0.
   $$
   Letting $\de\to 0$, we get that
   $$
   S(y,y)\le S(\s y,y)+A(y)-m_0
         \le S(\s y,y)+S(y,\s y)
         \le S(y,y).
   $$
   Thus $S(\s y,y)=-A(y)+m_0$ and $S(y,\s y)=A(y)-m_0$. Now
   \begin{align*}
   S(x,y)\ge S(x,\s y)+S(\s y,y)
         &=S(x,\s y)-A(y)+m_0
         \\
         &\ge S(x,y)+S(y,\s y)-A(y)-m_0
         \ge S(x,y).
  \end{align*}
  Hence $S(x,y)-S(x,\s y)=-A(y)+m_0$, and then
  $B(y)=A(y)-m_0+S(x,y)-S(x,\s y)=0$.

  \medskip  
  
  Now we prove that if $\mu$ is a maximizing measure then
  $\supp(\mu)\subseteq\fS$.  A proof of the following lemma
  is supplied below:

  \begin{Lemma} (Ma\~n\'e~\cite{Ma2})\label{L:G.5}
  
  Let $(X,{\mathcal B},\mu,f)$ be an ergodic measure preserving
  dynamical system and $F:X\to\re$ an integrable function.
  Given $A\in{\mathcal B}$ with $\mu(A)>0$, denote by $\hA$
  the set of points $x\in A$ such that for all $\e>0$ there 
  exists an integer $N>0$ such that $f^N(x)\in A$ and
  $$
  \Big|\zum_{k=0}^{N-1} F(f^k(x))-N\displaystyle\int F\,d\mu\Big|<\e.
  $$
  Then $\mu(A)=\mu(\hA)$. 
  \end{Lemma}             

  \medskip
  
  Let $\mu$ be a maximizing measure and $y\in\supp(\mu)$.
  Let $\de>0$, $z\in\SS$ and $n>0$ such that $d(y,z)<\de$,
  $d(\s^nz,y)<\de$ and
  $$
  \zum_{k=0}^{n-1}\b[A(\s^kz)-m_0\b]>-\de.
  $$
  The set $\{y,\s z,\s^2 z,\ldots,\s^{n-1}z,y\}$ is a
  $\de$-pseudo-orbit with 2 jumps. By lemma~\ref{L:G.2}, 
  there is $w\in\SS$ with $d(w,y)<\frac{\de}{1-\la}$, 
  $\s^n w= y$ and 
  $$
  \Big|\zum_{k=0}^{n-1}\b[A(\s^k w)-m_0\b]
  -\zum_{k=0}^{n-1}\b[A(\s^k z)-M_0\b]\Big|
  \le 2\, K\, \de^\a.
  $$  
  Hence
  $$
  S_{\frac{\de}{1-\la}}(y,y)\ge -\de - 2\, K\, \de^\a.
  $$
  Letting $\de\to 0$ we get that $S(y,y)=0$.
  
  \qed
  
  \bigskip

  \noindent{\bf Proof of Lemma~\ref{L:G.5}:}
  
  We may assume that $\int F\, d\mu = 0$. For $\e>0$ let
  $$
  A(\e):=\big\{\,p\in A\;\big\vert\;\exists N>0,\; f^N(p)\in A,\;
  \big|\textstyle\sum_{k=0}^{N-1} F(f^kp)\big|<\e\,\big\}.
  $$
  Let $x\in A$ be a point such that Birkhoff's Theorem holds for
  $F$ and the characteristic functions $1_A$ and $1_{A(\e)}$.
  It is enough to prove that $\mu(A(\e))=\mu(A)$ because
  $\hA=\cap_{n>0} A(1/n)$.
  
  Let $N_1<N_2<\cdots$ be the integers for which $F^{N_i}(x)\in A$.
 Define $\de(k)\ge 0$ by $N_k\cdot \de(k) =|\sum_{i=0}^{N_k-1}
 F(f^ix)|$. Then $\lim_{k\to+\infty}\de(k)=0$.
 
 Let $c_j:= \sum_{i=0}^{N_j-1} F(f^ix)$ and
 $$
 S(k):=\bigl\{\, 1\le j\le k-1\,\vert\, \not\exists \ell>j
 \text{ with }\abs{c_\ell - c_j}<\e\,\bigr\}.
 $$
 Then $\e\;\#S(k)\le 2\,\de(k)\, N_k$.
 
 If $j\notin S(k)$ then 
 $|c_\ell-c_j|=|\sum_{N_j}^{N_\ell -1} F(f^ix)|<\e$ for some
 $\ell>j$, hence $f^{N_j}(x)\in A(\e)$.
 
 We have that
 \begin{align*}
 \frac 1{N_k}\;\#\bigl\{\,0\le j <N_k\,\big\vert\,
 f^j(x)\in &A\setminus A(\e)\,\bigr\}
   \le \frac 1{N_k}\;\#S(k) \\
   &\le\frac 1{N_k}\cdot\frac{2\de(k)}{\e}\, N_k
   =\frac{2\,\de(k)}{\e}\overset{k}\longrightarrow 0.
 \end{align*}
 The choice of $x$ implies that $\mu(A\setminus A(\e))=0$.  
 
 \qed

  \bigskip
  \bigskip
  \bigskip

  To give an idea of how discontinuous the functions
  $S(x,y)$ and $y\mapsto S(x,y)$ ($x\notin\fS$) may be,
  we show the following proposition:
  
  \bigskip
  
  \begin{Proposition}\label{P:discontinuous}\quad
  
  Given $x\in\SS$ and $0<N\le\min\{k>0\,\vert\,\s^k(x)=x\,\}\le
  +\infty$, then
  $$
  S(x,\s^N x) = \zum_{k=0}^{N-1}\bigl[A(\s^k x)-m_0\bigr]
  $$
  and $S(x,x)=S(x,\s^N x)+S(\s^N x,x)$.
  \end{Proposition}
  
%   \bigskip
%   
%   \begin{Lemma}\label{2defs}\quad
%   The definition of $S_\de(x,y)$ may be replaced by
%   $$
%   \cS_\de(x,y) = \sup\Big\{\zum_{k=0}^{N-1}\bigl[A(\s^k z)-m_0\bigr]
%   \,\Big\vert\, d(\s^n z,y)<\de,\; d(x,z)<\de\,\Big\},
%   $$
%   giving the same limit $S(x,y)=\lim_{\de\to 0}\cS_\de(x,y)$.
%   \end{Lemma}
%   
%   \medskip
%   
%   \noindent{\bf Proof:}
%   If $d\s^n z,y)<\de$ and $d(x,z)<\de$, let $\phi_n:B(y,\e_0)\to\SS$ be
%   the branch of the inverse of $\s^n$ such that $\phi_n(\s^n z)=z$. Let
%   $w=\phi_n(y)$. Then $d(\s^kz,\s^kw)<\de\,\la^{n-k}$, and hence
%   $$
%   \zum_{k=0}^{n-1}[A(\s^k w)-m_0]
%   \ge \zum_{k=0}^{n-1}[A(\s^k z- m_0]
%   -\zum_{k=0}^{n-1}\Hold_\a(A)\,\de^a\,\la^{\a k}.
%   $$
%   Thus
%   $$
%   \cS_\de(x,y)\ge   S_\de(x,y)
%   \ge\cS_\de(x,y)-\frac{\Hold_\a(A)\,\de^a}{1-\la^a}.
%   $$
%   \qed
%   
  
  \medskip
  
  \noindent{\bf Proof:}
  
  Fix $x\in\SS$ and $N>0$ as in the statement of
  proposition~\ref{P:discontinuous}. Let $\e>0$ be small and
  $0<\de<\e$ such that if $d(z,x)<\de$ then
  \begin{equation}\label{E:d.1}
  d(\s^k z, \s^k x) < \e \quad\text{ for all } 0\le k \le N.
  \end{equation}
  Let $w\in\SS$ and $M>0$ be such that $d(w,x)<\de$, 
  $\s^M w=x$ and
  $$
  \zum_{k=0}^{M-1}\bigl[ A(\s^k w) - m_0\bigr] 
  \ge S_\de(x,x)-\de.
  $$
  If $0<2\e<\min\{\, d(\s^i x,\s^j x)\,\vert\, 0\le i<j\le N\, \}=:D$,
  then $M> N$ because for $0<k\le N$ we have that 
  $$
  d(\s^k w,\ x) \ge d(\s^k x,  x) -d(\s^k x,\s^k w)
                   > D-\e > \de.
  $$
  From~\eqref{E:d.1}, we have that 
  $$
  \zum_{k=0}^{N-1}\bigl[A(\s^k x)-m_0\bigr]
  \ge \zum_{k=0}^{N-1}\bigl[ A(\s^k w)-m_0\bigr]
      - N\,K\,\e^\a,
  $$
  where $K$ is an $\a$-H\"older constant for $A$. Then
  \begin{align*}
  S(x,\s^N x)&+S_\e(\s^N x,x)\ge 
    \zum_{k=0}^{N-1}\bigl[ A(\s^k x)-m_0\bigr]+S_\e(\s^Nx,x) 
    \\
    &\ge\zum_{k=0}^{N-1}\bigl[A(\s^k w)-m_0\bigr]-N\,K\,\e^\a 
           +\zum_{k=N}^{M-1}\bigl[A(\s^k w)-m_0\bigr]  
    \\
    &\ge S_\de(x,x)-\de-N\,K\,\e^\a.
  \end{align*}
  Letting $\e\to 0$, we have that
  \begin{align*}
  S(x,x)&\ge S(x,\s^N x)+S(\s^N x,x) \\
        &\ge\zum_{k=0}^{N-1}\bigl[A(\s^k x)-m_0\bigr]+S(\s^N x,x)\\
        &\ge S(x,x).
  \end{align*}
  \qed

  \bigskip
  \bigskip

  \refstepcounter{section}\label{CVS}
  \noindent{\large\bf \thesection. The continuously varying support property.}
  
  \bigskip
  \bigskip

    \noindent{\bf Definition:}
    We say that a pair $(A,\mu)\in\Ca\times\cM(\sigma)$
    has the {\it semi-continuously varying support property\/}
    if for any neighbourhood $U\subseteq\SS$ of $\supp(\mu)$
    there exists a neighbourhood $\cV\ni A$ of $A$ in the
    $C^0$-topology, such that if $\phi\in\cV$,
    and $\nu$ is a maximizing measure for $A+\phi$, then
    $\supp(\nu)\subseteq\ U$.
    
    \bigskip
    
    \begin{Lemma}\label{cvs+uniqueness}\quad
    
    If a function $A\in\Ca$ has a unique minimizing
    measure $\mu$ and the semi-continuously varying support
    property, then $\supp(\mu)$ is uniquely ergodic and
    $\mu$ has the continuously varying support property.
    \end{Lemma}
    
    \medskip
    
    \noindent{\bf Proof:}
    
    The unique ergodicity follows from item~[4] of
    proposition~\ref{P:G.3}. To obtain the continuously varying
    support property we have to show that the map $\Ca\ni A\to 
    \cM(\SS)$ is continuous in the strong topology. 
    By the hypothesis of semi-continuity, it is enough to prove that if
    $\psi_n\in C^\a(\SS,\re)$ and $\norm{\psi}_0\to 0$, then
    $\nu_n\to \mu$ weakly*, where $\nu_n$ is a maximizing measure for
    $A+\psi_n$. 
    
    Choose a limit $\tnu$ of a subsequence of $\nu_n$.
    Then $\int (L+\psi_n)\, d\mu \le \int
    (L+\psi_n)\, d\nu_n$ and hence $\int L\,d\mu\le\int L\, d\tnu$.
    Thus $\tnu$ is maximizing for $A$ and hence $\tnu=\mu$.    
    
    \qed
    
    \bigskip
    
    Theorem A combined with the following proposition
    give a proof of theorem C.
    
    \bigskip

    \begin{Proposition}\label{P:cvs}\quad
    
    Let $A^*\in\Ca$ admiting a unique maximizing measure
    $\mu^*$. Let $\Psi:\SS\to\re$ be a continuous function 
    such that $\Psi(x)=0$ for $x\in\supp(\mu^*)$ and $\Psi(x)<0$ for
    $x\notin\supp(\mu^*)$. Then $(A^*+\Psi,\mu^*)$ has the 
    semi-continuously varying support property.
    
    \end{Proposition}
    
    \medskip

  \noindent{\bf Proof:}
    
    Write $A:=A^*+\Psi$. 
    By lemma~\ref{cvs+uniqueness}, it is enoungh to prove the
    semi-continuosly varying support property.  
     Suppose that it does not hold.
    Then there is a neighbourhood $U$ of $\supp(\mu^*)$ and a sequence  
    $\langle A_n\rangle_{n\ge 0}\subset\Ca$ of H\"older functions 
    converging to $A$ and maximizing measures $\mu_n$ for $A_n$
    such that $K_n=\supp(\mu_n)\not\subseteq U$. We may assume that
    $\mu_n$ converges weakly to $\mu_\infty$ and $\langle
    K_n\rangle_{n\ge 0}$ converges in the Hausdorff metric to a
    compact set $K_\infty$.
    
    \medskip
    
    \noindent{\bf Step one:}
    Let $\la_n=\int A_n\; d\mu_n$ and $\la^*=\int A\; d\mu^*$. We
    prove that $\la_n\to\la^*$ and that $\langle\mu_n\rangle_{n\ge
    0}$ converges weakly* to $\mu^*$.
    
    We have that $\la^n\ge\int A_n\, d\mu^*$, hence
    $$
    \lim\inf_n \la_n\ge\int A\; d\mu^* = \la^*.
    $$
    Moreover,
    \begin{align*}
    \la^*=\int A\; d\mu^*
         &\ge \int A\; d\mu_n \\
         &\ge \int A_n\; d\mu_n -\norm{A-A_n}_0
         =\la_n - \norm{A-A_n}_0.
    \end{align*}
    Letting $n\to \infty$, we get that
    $\lim_n\sup \la_n \le\la^*$.        
    
 \medskip
 
 \noindent{\bf Step two:}
 We how that we can extend the coboundary 
 equation for $A$ to  $K_\infty$.
 
 Fix $\ov{x}\in\SS$. Let $V_n\in\Can$ be a function given by
 proposition~\ref{P:G.3}[5] for $A_n$. By adding a constant
 we may assume that $V_n(\ov{x})=0$. By proposition~\ref{P:G.3}[5],
 $\Hold_\a(V_n)$ is uniformly bounded on $n$. By Arzela-Ascoli
 theorem, there is a convergent subsequence
 $V_n\overset{\norm{\;}_0}{\longrightarrow}W$ to an $\a$-H\"older
 function $W$. Since
 $$
 A_n =\la_n + V_n - V_n\circ\s \qquad\text{ on }K_n, 
 $$  
 then
 \begin{equation} 
 A =\la^* + W - W\circ\s \qquad\text{ on }K_\infty. \label{cvs.1} 
 \end{equation}
 Similarly,
 \begin{equation}
 A \le\la^* + W - W\circ\s \qquad\text{ on all }\SS. \label{cvs.2}
 \end{equation}
 Since $\int A\,d\mu^*=\int A^*\, d\mu^*=\la^*$, then
 \begin{equation}\label{cvs.3}
 A=\la^*+W-W\circ\s \qquad\text{ on }\supp(\mu^*).
 \end{equation}
 
 \medskip
 
 \newpage
 \noindent{\bf Step three:}
 
 Since $K_n=\supp\mu_n$ then for all $x\in K_n$ there is a complete
 foward orbit in $K_n$ containing $x$, i.e. there is $\langle
 x_k\rangle_{k\in\Z}$ such that $x_0=x$ and $\s(x_k)=x_{k+1}$ for all
 $k\in\Z$. Then $K_\infty$ has also this property.
 
 Let $y\in K_\infty\setminus U$ and $\langle
 y_k\rangle_{k\in\Z}\subseteq K_\infty$ such that $y_0=y$ and
 $\s(y_k)=y_{k+1}$, $\forall k\in\Z$. By the cohomology
 property~\eqref{cvs.1}, any invariant measure supported on $K_\infty$
 is maximizing and thus it is $\mu^*$. Hence there are sequences
 $M,\;N\to +\infty$ such that
 $$
 \tfrac 1N\zum_{k=0}^{N-1}\de_{y_k}\overset{w^*}\longrightarrow\mu^*
 \qquad \text{ and }\qquad
 \tfrac 1M\zum_{k=-M}^{-1}\de_{y_k}\overset{w^*}\longrightarrow\mu^*,
 $$
 where $\de_y$ is the Dirac probability supported on $\{y\}$ and the
 convergences are in the weak* topology. In particular, we may assume 
 that $d(y_N,\supp\,\mu^*)\to 0$ and $d(y_{-M},\supp\,\mu^*)\to 0$.
 Since $\mu^*$ is uniquely minimizng, then in is ergodic.
 By the ergodicity of $\mu^*$, there is $z=z(N,M)\in\supp(\mu^*)$ and
 $K=K(N,M)>0$ such that $d(z,y_N)\to 0$ and $d(\s^Kz,y_{-M})\to 0$. The
 sequence $y_{-M},\ldots,y_0,\ldots,y_{N-1},z,\ldots,\s^{K-1}z$ is a
 closed $\e$-pseudo orbit with 2 jumps and with $e=\e_{N,M}\to 0$.

 Let $B=A-\la^*+W-W\circ\s\le 0$. By~\eqref{cvs.1} and~\eqref{cvs.3},
 $B=0$ on $K_\infty\cup\supp(\mu^*)$. By lemma~\ref{L:G.2}[2], there is a
 periodic point $p\in\SS$ such that $d(p,y_0)<\frac{\e}{1-\la^a}$ and
 \begin{align*}
 -\zum_{k=0}^{M+N+K-1}B(\s^kp)
 &=\zum_{k=-M}^{N-1}B(y_k)_+\zum_{k=0}^{K-1}B(\s^kz)
 -\zum_{k=0}^{M+N+K-1}B(\s^kp)  \\
 &<2\, K_1\,\e^\a.
 \end{align*}
 Now,
 \begin{align*}
 \zum_{k=0}^{M+N+K-1}\bigl[A^*(\s^kp)-\la^*\bigr]
 &=\zum_{k=0}^{M+N+K-1}B(\s^kp)
 -\zum_{k=0}^{M+N+K-1}\Psi(\s^kp)   \\
 &\ge -\Psi(p) + 2\, K_1\,\e^\a.
 \end{align*}
 Since $p\to y_0$ and $\Psi(y_0)<0$ then, for $\e>0$ small, we have
 that
 \begin{equation}\label{cvs.4}
  \zum_{k=0}^{M+N+K-1}\bigl[A^*(\s^kp)-\la^*\bigr]>0.
 \end{equation}
 If $\nu_p$ is the invariant measure supported on the positive orbit
 of $p$, then~\eqref{cvs.4} implies that $\int A^*\,d\nu_p>\la^*$.
 This contradicts the choice of $\mu^*$. 
 
 \qed

  \bigskip
  
  \begin{Remark}\label{R:cvs} 
  If in proposition~\ref{P:cvs} we need $B$ and $B+\Psi$ to have 
  pressure zero, we can replace $B+\Psi$ by $t\,(B+\Psi)$ such that
  $P\bigl(t\,(B+\Psi)\bigr)=0$. Since the function
  $f(t,\Psi)=P\bigl(t\,(B+\Psi)\bigr)$ is analytic on $\re\times
  C^a(\SS,\re)$, then $\Psi$ can be chosen $C^\a$-arbitrarily close to~$0$ 
  and $t$ arbitrarily close to $1$. In particular, $t\,(B+\Psi)$ can
  be made $C^\a$ arbitrarily close to $B$ for any $0<\a\le 1$.
  \end{Remark} 

  \bigskip
  \bigskip

  \refstepcounter{section}\label{PERIODIC}
  \noindent{\large \bf \thesection. Maximizing measures 
                                for generic potentials.}
  
  \bigskip
  \bigskip
  
  Let $\Ca$ be the set of $\a$-H\"older continuous functions
  $A:\SS\to\re$ which have topological entropy $P(A)=0$, endowed
  with the $\a$-H\"older norm $\lV A\rV_\a:=\norm{A}_0 + \norm{A}_\a$.
  Let $\Cap$ be the closure in the $\a$-H\"older topology of
  $\cup_{\ga>\a}C^\ga_0(\SS,\re)$.
  
  If $p\in\Fix\,\s^N$, let $\nu_p$ be the probabiliy measure defined
  by
  $$\int f\, d\nu_p = \frac 1N \, \zum_{k=0}^{N-1} f(\s^k p).
  $$
  For convenience of the reader we rephrase theorem B.

  \bigskip
  \bigskip

  \noindent{\bf Theorem B.}
  
  {\it
             Let $\cG_2\subset\Cap$ be the set of $A\in\Cap$ such that
             there is a neighbourhood $\cU\ni A$ such that
             for all $B\in\cU$, the unique maximizing measure 
             for $B$ is $\nu_p$. Then $\cG_2$ is open and dense in
             $\Cap$.
  }

  \bigskip

  \bigskip
  
  \noindent{\bf Proof:}
  
  Let $\cH\subset\Cap$ be the set of $A\in\Cap$ such that there is a
  unique maximizing measure for $A$ which is supported on a periodic
  orbit. By proposition~\ref{P:density}, the set $\cH$ is dense on
  $\Cap$. By proposition~\ref{P:cvs} and remark~\ref{R:cvs}
  there is a dense subset $\cA\subseteq\cH$ such that any 
  $A\in\cA$ has the semi-continuously varying support property. 
  Then $\cA$ is dense in $\Cap$. We show now that $\cA=\cG_2$ and, 
  in particular, that it is open on $\Cap$. Let $A\in\cA$ and 
  let $p\in\SS$ be a periodic
  point such that the maximizing measure for $A$ is $\nu_p$.
  There exists a neighbourhood $U$ of $\cO(p)$ such that the
  unique invariant measure supported on $U$ is $\nu_p$. Since $A$ has 
  the continuously varying support property, then there is a
  neighbourhood $\cU(A)\subset\Cap$ such that the (unique) maximizing 
  measure for any $B\in\cU(A)$ is $\nu_p$.

  \qed

  \bigskip

  \begin{Proposition}\label{P:density}\quad
  
  The set $\cH$ of functions $A\in\Cap$ such that 
  $A$ has a unique minimizing measure and this measure 
  is supported on a periodic orbit is dense
  on $\Cap$.
   
  \end{Proposition}
  
  \bigskip
  
  \noindent{\bf Proof:}

  Let $F\in\Cap$, then for any $\rho>0$ there is $\a<\ga<1$ and
  $A\in\Cga$ such that $\lV A-F\rV_\a<\rho$. Let $\a<\be<\ga$,
  we will find $G=A+\Psi\in\Cb$ such that  $\lV \Psi\rV_\be<\rho$.
  Then $G\in\Cap$ and $\lV \Psi\rV_\a\le\lV \Psi\rV_\be<\rho$. 
  
  Let $\mu$ be a maximizing measure for $A$.  
 Suppose that there are no periodic points on $\supp(\mu)$. Then
 for all $n>0$, $\min_{z\in\supp(\mu)}d(z,\s^n z)>0$. Because otherwise
 $\Fix\,\s^n\cap\supp(\mu)\ne\0$. We will first find a periodic point
 sufficiently close to $\supp(\mu)$.
 
 Let $\eta:=\frac 12(1-\la)$ and let $K>0$ be such that
 \begin{align}
 \frac 1{1-\la^K}<\tfrac 32, \label{E:f.1}\\
 1-\frac{\la+\la^K}{1-\la^K}>\eta. \label{E:f.2}
 \end{align}
 and let $D>0$ be such that
 $$
 \min\{\, d(z,\s^j z)\,\vert\, z\in\supp(\mu),\; 0<j\le K\,\}>3D.
 $$
 Since  $\cup_{n\le K}\Fix\,\s^n$ is finite, there is $0<\e_1<D$ such
 that 
 \begin{equation}\label{E:f.3}
 \inf\{\, d(z,\s^jz)\,\vert\, d(z,\supp(\mu))<2\e_1,\; 0<j\le
 K\,\}>2D.
 \end{equation}
 Given $0<\e<\e_1$, let $z\in\supp(\mu)$ and $n>0$ be such that
 $$
 d:=d(z,\s^n z)=\min\{\, d(\s^iz,\s^jz)\,\vert\, 0\le i<j\le n\,\}<\e.
 $$
 By~\eqref{E:f.3}, we have that $n>K$. Using lemma~\ref{L:G.2}, we get
 that there exists $p\in\Fix\,\s^n$ such that $d(p,\s^nz)\le \frac
 d{1-\la^n}$ and for $0\le j\le n$,
 \begin{equation}\label{E:f.4}
 d(\s^jp,\s^jz)\le \frac{d\,\la^{n-j}}{1-\la^n}
               \le\tfrac 32 \, d < \tfrac 32 \,\e< 2\, \e_1.
 \end{equation}               
  Given $0\le i<j\le n-1$ by~\eqref{E:f.4} and~\eqref{E:f.3}, we have
  that
  \begin{equation}\label{E:f.5}
  d(\s^ip,\s^jp)>2D>\eta\, d \quad\text{ if } j\le i+K,
  \end{equation}
  and using~\eqref{E:f.4} and~\eqref{E:f.2},
  \begin{align}
  d(\s^ip,\s^jp) 
    &\ge d(\s^iz,\s^jz)- d(\s^iz,\s^ip)-d(\s^jz,\s^jp) \notag\\
    &> d -
    \frac{\la^{n-i}\,d}{1-\la^n}-\frac{\la^{n-j}}{1-\la^n}\notag\\
    &>\left[1-\frac{\la+\la^K}{1-\la^n}\right]\, d > \eta\, d
     \quad\text{ if }i+K\le j.\label{E:f.6}
  \end{align}

  Fix $q\in\supp(\mu)\subseteq\fS$ and $W:\SS\to\re$, 
  $W(y)=S(q,y)$. Then $W$ is
  $\ga$-H\"older continuous and
  $$
  W (\s x) - W(x) \ge S(x,\s x)\ge A - m_0 \quad\text{ for all }
  x\in\SS.
  $$
  Hence
  $$
  W \circ \s - W = A - m_0 \qquad \text{ on } \supp(\mu).
  $$
  Let $B(x):=A(x)-m_0+W(x)-W(\s x)\le 0$. Let $K_1>0$ be an
  $\ga$-H\"older constant for $B$. Let 
  $$
  \de = \tfrac 14 \,\eta\, d \quad\text{ and }\quad
  Q= K_1\left[\tfrac4\eta\right]^\ga > K_1.
  $$
  If $d(x,y)<\de$ and $0<\be<a$, 
  then
  $$
  \lv B(x)-B(y)\rv \le K_1 \,d(x,y)^\ga
                   < K_1\,\de^{\ga-\be}\, d(x,y)^\be.
  $$
  For $x\in\SS$, define
  $\abs{x}:=\min\{d(x,\s^kp)\,\vert\, k=1,\ldots,n\}$
  and $p_x=\s^kp$ such that $d(x,p_x)=\abs{x}$.
  Let
  \begin{equation}\label{E:f.phi}
  \Phi(x) = \max\Big\{\, 0\;,\; \left[3\,Q\, \de^{\ga-\be}
  -\tfrac{B(p_x)}{\de^\be}\right](\de^\be-\abs{x}^\be)\Big\}.
  \end{equation}
  We show that $\max_x B(x)+\Phi(x)= Q\de^\ga =
  B(\s^kx)+\Phi(\s^kx)$ for all $k=1,\ldots,n$.
  Indeed, for $\abs{x}<\de$
  $$
  \lv B(x)-B(p_x)\rv\le K_1\,\abs{x}^\ga
  \le Q\, \abs{x}^\ga
  \le (Q\,\de^{\ga-\be})\,\abs{x}^\be.
  $$
  If $p_x=\s^ip$, then
  \begin{align*}
  \lv B(p_x)\rv 
    &\le \lv B(\s^iz)\rv + K_1\, d(\s^iz,p_x)^\ga \\
    &\le 0+ K_1\, \left[\tfrac\la{1-\la^n}\right]^\ga\, d^\ga
     \le K_1\, d^\ga \le Q\, \de^\ga.
  \end{align*}

  Hence
  \begin{align*}
  B(x)+\Phi(x) 
  &\begin{aligned}[t]
   \le B(p_x) + Q\,\de^{\ga-\be}\,\abs{x}^\be
             &+ 3\,Q\de^\ga - B(p_x) 
             \\
             &- 3\, Q\, \de^{\ga-\be}\,\abs{x}^\be
             + \frac{Q\,\de^\ga}{\de^\be}\, \abs{x}^\be
             \end{aligned}
  \\
  &\le 3\,Q\,\de^\ga - Q \de^{\ga-\be}\,\abs{x}^\be
   \le 3\, Q\,\de^\ga.
   \end{align*}
  Also $B(p_x)+\Phi(p_x)=3\, Q\, \de^a$ and $B(x)+\Phi(x)=B(x)\le 0< 3\,
  Q\,\de^a$ for $\abs{x}>\de$.
  
  If $\nu\ne\nu_p$ is a $\s$-invariant probability, we have
  that
  $$
  \int A\;d\nu = \int\b[A+W - W\circ\s\b]\;d\nu
  \le \int B(x)\;d\nu_p +m_0
  <\int B\; d\nu_p  +m_0.
  $$
  
  We now prove that the $\be$-H\"older norm of $\Phi$ can be made
  arbitrarily small. We have that
  $$
  \lV \Phi\rV_0 := \sup_{x\in\SS}\lv\Phi(x)\rv
     \le 3\,Q\,\de^\ga + \max_{0\le i\le n-1}\lv B(\s^ip)\rv
     \le 4\,Q\,\de^\ga.
  $$   
  Observe that if $d(x,y)\le \de$ and $\abs{y}\le\de$
  then by~\eqref{E:f.5}
  and~\eqref{E:f.6} we have that $p_x=p_y$. If
  $\abs{y}\le\abs{x}\le 2\de$ and $0<\be<1$, then
  $$
  \abs{x}^\be-\abs{y}^\be
     \le \bigl(\abs{x}-\abs{y}\bigr)^\be
      \le d(x,y)^\be\,.  
  $$
  And if $\abs{y}\le\abs{x}\le\de$, then
  \begin{align*}
  \abs{\Phi(x)-\Phi(y)}
     &\le (3\,Q\,\de^{\ga-\be}+Q\,\de^{\ga-\be})
          \bigl(\abs{x}^\be-\abs{y}^\be\bigr)   
     \\
     &\le 4\, Q\,\de^{\ga-\be}\,d(x,y)^\be.
  \end{align*}        
  If $\abs{y}\le\de<\abs{x}$ and $d(x,y)\le \de$, then
  \begin{align*}
  \abs{\Phi(x)-\Phi(y)} 
  &\le 4\,Q\,\de^{\ga-\be}\,\bigl(\de^\be-\abs{y}^\be\bigr)
  \le 4\,Q\,\de^{\ga-\be}\,\bigl(\abs{x}^\be-\abs{y}^\be\bigr)
  \\
  &\le 4\, Q\,\de^{\ga-\be}\,d(x,y)^\be.
  \end{align*}
  If $d(x,y)\ge \de$ then
  $$
  \abs{\Phi(x)-\Phi(y)}
    \le \abs{\Phi(x)}\le 4\, Q\, \de 
    \le 4\, Q\, \de^{\ga-\be}\, d(x,y)^\be.
  $$  
  Hence
  $$
  \Hold_\ga(\Phi):=\sup_{0<d(x,y)\le 1}
  \frac{\abs{\Phi(x)-\Phi(y)}}{d(x,y)^\be}
  \le 4\, Q\,\de^{\ga-\be}.
  $$
  If we let $\e\to 0$ then $\de\to 0$, $\norm{\Phi}_0\to 0$ and
  $\Hold_\be(\Phi)\to 0$ for any $0<\be<\min\{1,\ga\}$.

  In the case when there is a periodic point
  $p\in\Fix\,\s^n\cap\supp(\mu)\ne\0$, choose
  $D>0$ such that $d(\s^ip,\s^jp)>2D$ for $0\le i<j\le n-1$ and define
  $\Phi(x)$ by the same formula as~\eqref{E:f.phi}. In this case
  $B(p_x)\equiv 0$. The rest of the proof is the same.

  Finally, we need to pertub $A$ among the $\be$-H\"older functions 
  with pressure zero. Let $t=t(\Phi)\in\re$ be such that
  $P(A+\Phi+t(\Phi))=0$. Since $\abs{P(A+\Phi)-P(A)}\le\norm{\Phi}_0$,
  $P(A+\Phi+t)=P(A+\Phi)+t$ and $P(A)=0$, then
  $\abs{t}\le\norm{\Phi}_0$. The perturbing function
  $\Psi=\Phi+t(\Phi)$ has $\norm{\Psi}_0\le 2\,\norm{\Phi}_0\le 4\,
  Q\,\de^\ga$ and $\Hold_\be(\Psi)=\Hold_\be(\Phi)\le 4\, Q\,\de^{\ga-\be}$.
  
  \qed

  \bigskip
  \bigskip

  \refstepcounter{section}\label{THERMO}
  \noindent{\large\bf \thesection. Relations with the Thermodynamic
                   Formalism.}
  \bigskip
  \bigskip

 \noindent{\bf Proof of Proposition D.}
 
 Let $\hmu_t$ be the equilibrium state for $tA$. Suppose that
 $\hmu_t$ does not converges weakly* to $\mu_A$, then for $\e>0$ small
 and a
 subsequence $t_n$,
 $$
 0< \int A\; d\hmu_{t_n} < \int A\; d\mu_A -\e.
 $$
 Take $t_n$ large enough such that $t_n\e-h_{top}(\s)> 0$. Then
 \begin{align*}
 h(\mu_A)+t_n\int A\;d\mu_A
 &\ge h(\mu_A) - t_n\Big[\int A\, d\hmu_{t_n}+\e\Big] \\
 &\ge h(\hmu_{t_n})-h_{top}(\s) +t_n\,\e +t_n\int A\, d\hmu_{t_n} \\
 &> h(\hmu_{t_n}) +t_n\int A\, d\hmu_{t_n}.
 \end{align*}
 \qed

  \bigskip
  \bigskip
  
  %\newpage
  \refstepcounter{section}\label{CIRCLE}
  \noindent{\large\bf \thesection. Expanding maps of the circle.}
  
  \bigskip
  \bigskip

 In this section we prove theorems A1 and B1. The idea is to
 show in proposition~\ref{homeo} below, a homeomorphism
 among the $C^{1+\a}$ expanding dynamics on $S^1$ and
 $C^\a$ functions on the correspondig shift $\SS$, and then
 to apply theroems~A and~B.
  
 \medskip
 
 Consider a point $y_0\in S^1$. In order to prove  theorems~A1 and B1, 
 it is enough to prove their claims for the class of
 maps $f\in\cF(\a)$ (resp. $\cF(\a+)$) 
 that fix the point $y_0$. We will also denote by $\cF(\a)$ (resp.
 $\cF(\a+)$) this new class of maps.
 
 We need to consider an abstract model that will be played by the
 transformation $T:S^1\to S^1$, given by $T(x)= 2x$~(mod 1). This map
 is equivalent to the full one-sided shift in two symbols 
 with identifications. We will use the diadic notation 
 for points in the circle without stressing the equivalence 
 of both systems.
 
 We call $x_0$ the fixed point of $T$. Given $f$, we will define
 a bi-H\"older map $\tf$ which conjugates $f$ and $T$, that is,
 $f\circ\tf=\tf\circ T$. In particular $\tf(y_0)=x_0$. 

 \medskip

 {\bf Construction of $\tf$:} \quad 
 
 Given a map $f$, let $z$ be the
 unique pre-image of $y_0$ different from $y_0$. Each point $t\in
 S^1$, $t\ne y_0$, has two different preimages in
 $S^1\setminus\{y_0\}$. These preimages $t_0$ and $t_1$ aqre ordered
 by the order of the interval $S^1\setminus\{y_0\}$, that is,
 $t_0<t_1$.
 
 We will order and code all pre-images $z_{\a_1,\a_2,\ldots,\a_n}(f)$
 (where $\a_i\in\{0,1\}$ and $n\in\na$) of $z$ in the following way:
 if $z_{\a_1,\a_2,\ldots,\a_n}(f)$ is defined, then
 $z_{\a_1,\a_2,\ldots,\a_n,0}(f)$ and $z_{\a_1,\a_2,\ldots,\a_n,1}(f)$
 are ordered by the previous procedure.
 
 We do the same for $T$ (substituting $y_0$ by $x_0$) and obtain a set
 of coded points $z_{\a_1,\a_2,\ldots,\a_n}(T)$ where $\a_i\in\{0,1\}$
 and $n\in\na$. Denote by $Z(f)$ and $Z(T)$ the set of preimages
 defined above respectively for $f$ and $T$.

 Define first $\tf$ in these points, by associating the corresponding
 points $Z(f)$ and $Z(T)$ with the same code. Then $\tf$ extends
 continuously to $S^1$ in a unique way, because both sets of preimages
 are dense on $S^1$. The map $\tf$ is a homeomorphism. By usual
 bounded distortion arguments, we obtain that $\tf$ is bi-H\"older.
 
 Consider the set $\Hl(\a)$ of $\a$-H\"older continuous functions
 $A:S^1\to\re$ wich are smaller than $-\log\la$. Observe that
 the topological pressure of $-\log f'\circ\tf$ is zero (see~\cite{PP}
 for a definition of topological pressure). Denote by $\Hlo(\a)$ the
 set
 of functions in $\Hl(a)$ with topological pressure zero and let
 $\Hlo(\a+)$ be the clousure in the $C^\a$-topology of
 $\cup_{\be>\a}\Hlo(\a)$.
 
 Define the transformation $\cG:\cF(\a)\to\Hlo(\a)$ by $\cG(f)=-\log
 f'\circ\tf$, where $f\in\cF_\la$. Similarly define
 $\cG:\cF(\a+)\to\Hlo(a+)$. Observe that $\tf$ depends on $f$
 in the definition of $\cG$.
 
 Theorems A1 and B1 follow from theorems A, B and the next proposition:

 \bigskip
 
 \begin{Proposition}\label{homeo}
 The transformations $\cG$ are homeomorphisms from $\cF(\a)$ [resp.
 $\cF(\a+)$] (with the $C^{1+\a}$ distance) to $\Hlo$ [resp.
 $\Hlo(\a+)$]
 (with the $C^\a$ distance).
 \end{Proposition}
 
 \medskip
 
 \noindent{\bf Proof:}  

 We shall prove that $\cG:\cF(\a)\to\Hlo(\a)$ is a homeomorphism
 for any $0<\a<1$. This implies that $\cG:\cF(a+)\to\Hlo(\a+)$ is
 a homeoporhism for any $0<\a<1$.
 
 We show first that $\cG$ is surjective.
 We have to find $f$ and $\tf$ as above for each given $A\in\Hlo$.
 
 Denote by $K(T)$ the set of invariant measures for $T$.
 For a given H\"older potential $A$ with pressure zero, denote by
 $\hm_A$ the eigenmeasure of the dual of the Ruelle-Perron-Frobenius
 operator of the potential $A$, that is $\cL_A^*\hm_A=\hm_A$
 (see~\cite{PP} for references on Thermodinamic Formalism). Note that
 the maximal eigenvalue of $\cL_A^*$ is $1$, because the pressure of
 $A$ is zero and that $\hm_A$ is not necessarily an invariant
 measure  in $K(T)$.
 
 Now we define a H\"older homeomorphism $\ta:S^1\to S^1$.
 By definition $\ta(x_0):=y_0$. For $x\ne x_0$ define
 $\ta(x)=y$ in such way that $\text{length}(y_0,y)=\hm_A(x_0,x)$.
 The map $\ta$ is well defined because $\hm_A$ is a probability
 with no atoms which is positive on open sets and the circle is
 oriented and has lenght one.
 
 Let $f=\ta\circ T\circ\ta^{-1}$. Since $\ta$ preserves orientation,
 then the two sets of preimages $Z(f)$ and $Z(T)$ are ordered in the
 same way. This proves that $\ta=\tf$.
 
 The Jacobian of $T$ with respect to the measure $\hm_A$ is $\E^{-A}$.
 By definition, the pushed measure of $\hm_A$ by $\ta$ is the Lebesgue
 measure. Since $f$ was defines by the change of coordinates $\ta$,
 then $f'$, the Jacobian of $f$ satisfies $f'=\E^{-A}\,\theta_A^{-1}$.
 Therefore $f'$ exists and it is H\"older. This shows that $\cG$ is
 surjective.
 
 Now we show that $\cG$ is injective. Suppose that two maps $f$ and
 $g$ satisfy $\cG=A_f=A_g=\cG(g)$. Consider the respective changes of
 coordinates $\tf$ and $\tg$.
 
 Note that $h=\theta_f^{-1}\circ\tg$ conjugates $f$ and $g$, because
 $\tf$ conjugates $f$ and $T$ and $\tg$ conjugates $g$ and $T$. Since
 $\tf=\tg$, because $A_f=A_g$, then $g$ is the identity and hence
 $f=g$. This implies that $\cG$ is injective.
 
 From the definition of $\cG$ and the reasoning above, it is easy to
 see that the map $\cG$ is an homeomorphism.
 
 \qed


\bigskip
 
   
    
               
 \begin{thebibliography}{99}

 \bibitem{A} G. Atkinson. {\it Recurrence of cocycles and random
            walks.\/} J. London Math. Soc. {\bf (2)}, N.{\bf 13}
            (1976), 486-488.

 \bibitem{C} M.J. Dias-Carneiro. {\it On minimizing measures of the
 action of autononomous Lagrangians}. Nonlinearity~{\bf 8}, (1995),
 1077-1085.
 
 \bibitem{CDI}  G. Contreras, J. Delgado \& R. Iturriaga.
                {\it Lagrangian Flows: The dynamics of globally
                minimizing orbits - II.\/} Bol. Soc. Bras. Mat.
                Vol. {\bf 28}, N.2, 155-196, (1997). 
                Also available via internet at
                ``http://www.ma.utexas.edu/mp\_arc''.
                
 \bibitem{CIPP} G. Contreras, R. Iturriaga, G.P. Paternain, M.
                Paternain. {\it Lagrangian graphs, minimizing measures
                and Ma\~n\'e's critical values}. To appear GAFA.
                Also available via internet at
                ``http://www.ma.utexas.edu/mp\_arc''.
                               
 \bibitem{Fa1} A. Fathi. Th\'eor\`eme KAM faible et th\'eorie de Mather
               sur les syst\`emes lagrangiens. C. R. Acad. Sci. Paris,
               t. {\bf  324} S\'erie {\bf I}, 1043-1046, (1997).
               
 \bibitem{Fa2} A. Fathi. Th\'eor\`eme KAM faible et th\'eorie de
               Mather sur les syst\`emes lagrangiens II. Preprint
               UMPA, ENS-Lyon.              

 \bibitem{FM}  G. Forni \& J. Mather. {\it Action minimizing orbits in
                Hamiltonian systems}, in Lect. Notes in Math. 1589,
                Springer Verlag. (1994), 92-186.

 \bibitem{G}  C. Grillenberger. {\it Constructions of strictly ergodic
              systems I. Given entropy}. Z. Wahrscheinlichkeitstheorie
              verw Geb, {\bf 25}, (1973), 323-334.
                
 \bibitem{H}   B. Hunt. {\it Maximal local Lyapunov dimension bounds
               the box dimensionof chaotic attractors}.
               Nonlinearity~{\bf 9}, (1996), 845-852.             
                       
 \bibitem{L} A. Lopes. {\it Dimension spectra and a mathematical model
             for phase transitions}. Advances in Applied Math.
             vol.~{\bf 11}, N.{\bf 4}, (1990), 475-502.
             
 
 \bibitem{Ma1}    R. Ma\~n\'e. {\it The Hausdorff dimension of
                  horseshoes of diffeomorphismsof surfaces}. 
                  Bol. Soc. Bars. Mat., vol.{\bf 20}, N.{\bf 2},
                  (1990), 1-24.
                  
 
 \bibitem{Ma2}     R. Ma\~n\'e.
                   {\it Generic properties and problems of minimizing
                   measures of lagrangian systems.}
                   Nonlinearity, {\bf 9}, (1996), no.2, 273-310.
 
 \bibitem{Ma3}     R. Ma\~n\'e.
                   {\it Lagrangian Flows: The Dynamics of Globally
                   Minimizing Orbits.} In Int. Congress on Dyn. Sys.
                   in Montevideo (a tribute to Ricardo Ma\~n\'e),
                   F. Ledrappier, J.L. Lewowicz, S. Newhouse eds., 
                   Pitman Research Notes in Math. {\bf 362} (1996)
                   120-131. Reprinted in Bol. Soc. Bras. Mat. Vol~{\bf
                   28}, N. 2, 141-153.

 \bibitem{PP} W. Parry \& M. Pollicott.  Zeta function and the
                 periodic orbit structure of hyperbolic dynamics.
                 Asterisque vol.~{\bf 187-188}, (1990).
                 
 \bibitem{PUZ} F. Przytycki, M. Urbanski \& A. Zdunik,
               {\it Harmonic, Gibbs and Hausdorff measures on
               repellers for holomorphic maps}. Annals of Math.~{\bf
               130}, (1989), 1-14.
               
 \bibitem{R} R. Rockafellar.  Convex Analysis. Princeton Univ. Press.
             1970.              
 
 \bibitem{Te} R. Teman. Infinite-dimensional dynamical systems in
              mechanics and physics. Springer Verlag (1993).
              
 \bibitem{T} Ph. Thieullen. {\it Entropy and the Hausdorff dimension
             for infinite-dimensional dynamical systems}.  Journal of
             Dynamics and Differential Equations. Vol.{\bf 4}, N.{\bf
             1}, (1992), 127-159.           
 
      
 \end{thebibliography}
\end{document}